\documentclass[final, a4paper, 10pt]{article}

\usepackage{theorem}
\usepackage{amsfonts}
\usepackage{amsbsy}
\usepackage{amsmath}
\usepackage{latexsym}
\usepackage{verbatim}

\theoremstyle{plain}
\newtheorem{thm}{Theorem}[section]
\newtheorem{cor}[thm]{Corollary}
\newtheorem{lem}[thm]{Lemma}
\newtheorem{prop}[thm]{Proposition}
\newtheorem{rem}[thm]{Remark}
\newenvironment{prf}[1][Proof]{\textbf{#1.} } {$\square$ \bigskip}

\numberwithin{equation}{section}

\newcommand{\MTRX} {\mathbf{M}}                         %
\newcommand{\IC} {\mathbb{C}}                           %
\newcommand{\IR} {\mathbb{R}}                           %
\newcommand{\IT} {\mathbb{T}}                           %
\newcommand{\IZ} {\mathbb{Z}}                           %
\newcommand{\BDD}[1]{\mathbf B \big( #1 \big) }         %
\newcommand{\CPT}[1]{\mathbf K \big( #1 \big) }         %
\newcommand{\BM} {\mathbf M}                            %
\newcommand{\CA} {\mathcal{A}}                          %
\newcommand{\CC} {\mathcal{C}}                          %
\newcommand{\CH} {\mathcal{H}}                          %
\newcommand{\CK} {\mathcal{K}}                          %
\newcommand{\CL} {\mathcal{L}}                          %
\newcommand{\CM} {\mathcal{M}}                          %
\newcommand{\CP} {\mathcal{P}}                          %
\newcommand{\CT} {\mathcal{T}}                          %
\newcommand{\im} {\mathrm{i}}                           %
\newcommand{\trace} {\mathrm{trace}}                    %
\newcommand{\tr} {\mathrm{tr}}                          %
\newcommand{\Tr} {\mathrm{Tr}}                          %
\newcommand{\ind} {\mathrm{index}}                      %
\newcommand{\df} {\mathrm{d}}                           %
\newcommand{\Inv} {\mathrm{Inv}}                        %
\newcommand{\Com} {\mathrm{Com}}                        %
\newcommand{\cstr} {\mathrm{C}^*}                       %
\newcommand{\ot} {\otimes}                              %
\newcommand{\by}{{\times}}                              %
\newcommand{\al} {\alpha}                               %
\newcommand{\la} {\lambda}                              %
\newcommand{\ro} {\varrho}                              %
\newcommand{\vfi} {\varphi}                             %
\newcommand{\om} {\omega}                               %
\newcommand{\OM} {\Omega}                               %
\newcommand{\norm}[1]{\left\Vert#1\right\Vert}          %
\newcommand{\set}[1]{\left\{#1\right\}}                 %
\newcommand{\innerp}[2]{\left<#1 \mid #2\right>}        %
\begin{document}

\title{Non--Vanishing functions and Toeplitz Operators on Tube--Type Domains}

\author{Adel B. Badi}

%
\maketitle

%

\begin{abstract}
We prove an index theorem for Toeplitz operators on irreducible tube--type domains and we extend our results to Toeplitz operators with matrix symbols. In order to prove our index theorem, we proved a result asserting that a non--vanishing function on the Shilov boundary of a tube--type bounded symmetric domain, not necessarily irreducible, is equal to a unimodular function defined as the product of powers of generic norms times an exponential function.
\end{abstract}

%

\section[Introduction]{Introduction}
\label{sec0}

The well--known I. Gohberg and M. Krein index theorem for Toeplitz operators on the unit circle is a relation between the (analytic) Fredholm index of a Toeplitz operator and the winding number of its symbol (the topological index). The theory of Toeplitz operators have been extended and studied on many generalized Hardy spaces. However, in most cases there is no result similar to Gohberg--Krein index theorem.

G. J. Murphy has introduced an elegant index theory for generalized Toeplitz operators in \cite{MUR1} in attempt to generalize Gohberg--Krein theorem. Murphy's theory is based on $\cstr$--algebras where he used orthogonal projections to define Toeplitz operators and traces to define the topological and analytical indices. He showed that the classical Gohberg--Krein index theorem for Toeplitz operators on the unit circle can be derived from his results and he proved an index theorem for Toeplitz operators on the group of unitary $2 \by 2$ matrices. Moreover, he showed that the Gohberg--Krein theorem for Toeplitz operators with matrix symbols on the unit circle follows from his results too. He has proved an index theorem for Toeplitz operators with matrix symbols on the group of $2\by2$ unitary matrices. In \cite{ADBA1}, we used Murphy's results to prove an index theorem for Toeplitz operators on the quarter--plane and we extended our results to the related Toeplitz operators with matrix symbols.

In this paper we  prove an index theorem for Toeplitz operators on irreducible tube--type bounded symmetric domains. Our results can be viewed as an extension of Murphy's result for Toeplitz operators on the compact group of unitary $2\by2$ matrices. In fact, this group is the Shilov boundary of an irreducible tube--type bounded symmetric domain. We prove our index theorem using H. Upmeier's theory for Toeplitz operators on bounded symmetric domains; see \cite{UPM1}. All of the examples in \cite{MUR1} are of Toeplitz operators on connected compact groups, hence, the calculation of the topological index is done using the result of E. van Kampen in \cite{VKM1} which is an extension of a result of H. Bohr in \cite{BOH1}. In this paper the symbol algebra is defined on connected compact spaces which are not necessarily groups. Thus, we need a result on these spaces similar to van Kampen's result. This result asserts that a non--vanishing function on the Shilov boundary of a tube--type bounded symmetric domains is equal to a unimodular function defined as the product of powers of generic norms times an exponential function. Recall that a complex--valued function is called {\em unimodular} if its absolute value is $1$. This result is valid for all tube--type domains, the reducible and irreducible ones. We prove this result in its full generality in spite of the fact that our index theorem is proved only for Toeplitz operators on irreducible tube--type domains.

In Section \ref{sec1}, we introduce tube--type domains and some Jordan triple basic concepts like tripotents and generic norms which are necessary in our work. Bounded symmetric domains have very rich geometric and Jordan algebraic theories, however, only those properties needed in our work will be given. The structure of the Shilov boundary of a tube--type domain is studied and some of its important properties are given. We also give the topological properties of a compact subset of the Shilov boundary called the reduced Shilov domain.

Section \ref{sec2} is dedicated to the proof of our result on non--vanishing functions on Shilov boundaries of tube--type domains. We study Hardy spaces on irreducible tube--type domains in Section \ref{sec3} where we depended on the results obtained by Upmeier. Finally, we prove our index theorem for Toeplitz operators on irreducible tube--type domains in Section \ref{sec4} and we extend it to the related Toeplitz operators with matrix symbols.

%

\section[Jordan Triples, Tripotents and Generic Norms]{Jordan Triples, Tripotents and Generic Norms}
\label{sec1}

Let $Z$ be a finite--dimensional vector space on the field $\IC$ of complex numbers provided with the unique norm topology. Polynomials can be defined on $Z$ without reference to a certain basis and we can define holomorphic functions using polynomials \cite[Sec. 1.1]{UPM1}.

 Let $D$ be a domain in $Z$, that is, a connected non--empty open subset of $Z$. A mapping $g: D \to D$ is said to be {\em biholomorphic} if it is bijective and $g, g^{-1}$ are holomorphic. If $D$ is bounded and for every $z \in D$ there exists a biholomorphic mapping $g_z$ with $z$ as an isolated fixed point of $g_z$ and $g_z \circ g_z$ is the identity map on $D$, then $D$ is called a {\em bounded symmetric domain}.

Any bounded symmetric domain is biholomorphically equivalent to a unique (up to a linear isomorphism) circled bounded symmetric domain \cite[Theorem 1.6]{LOO1}. Recall that the set $D$ is {\em circled} if it is invariant under the map $Z \to Z, \ z \mapsto \la z$, $\la \in \IT$ and $0 \in D$, where $\IT$ is the unit circle in $\IC$. We will assume that all bounded symmetric domains are circled.

A {\em Jordan triple product} on $Z$ is a mapping $Z^3 \to Z$, $(x, y, z) \mapsto \set{xyz}$ which is linear in the first variable, conjugate--linear in the second variable, symmetric in the first and the third variables and
\[
\set{x y \set{z u v}} + \set{z \set{y x u} v} = \set{\set{x y z} u v} + \set{z u \set{x y v}},
\]
for all $x, y, z, u, v \in Z$. If $\set{z z z} = tz$ for $z \in Z$ implies $t > 0$ or $z = 0$, then we call $Z$ a {\em positive hermitian Jordan triple system}, PJT for short.

Let $Z$ be a PJT and define the mapping $L(x, y): Z \to Z, \ z \mapsto \set{x y z}$ for any $x, y \in Z$. Then the sesquilinear form
\begin{eqnarray}
\innerp{x}{y} = \trace \big( L(x, y) \big), \ \ x, y \in Z    \label{eqn5000}
\end{eqnarray}%
 is a positive--definite inner product
on $Z$. If we denote by $1$ the identity map on $Z$, then the set
\begin{eqnarray}
\set{z \in Z \mid 1 - L(z, z) \textrm{ is positive definite}} \label{eqn5001}
\end{eqnarray}%
is a bounded symmetric domain. Conversely, for every bounded symmetric domain $D \subseteq Z$, there exists a Jordan triple product on $Z$ such that $Z$ is a PJT and $D$ is equal to the set defined in (\ref{eqn5001}). Moreover, there is a unique Jordan triple product such that $D$ is the open unit ball with respect to the so--called spectral norm on $Z$. Henceforth, we will always relate this unique Jordan triple product on $Z$ to any given bounded symmetric domain $D \subseteq Z$.

Let $x$ be an element in the bounded symmetric domain $D \subseteq Z$. If the conjugate--linear mapping $Q_x: Z \to Z, \ z \mapsto \set{x z x}$ is invertible, then we call $x$ an {\em invertible} element. The domain $D$ is said to be of {\em tube--type} if $Z$ contains an invertible element. A bounded symmetric domain is said to be {\em irreducible} if it is not a direct product of two bounded symmetric domains. Let $D_i$ be a bounded symmetric domain with the related PJT $Z_i$, where $i = 1, 2$. By \cite[4.11]{LOO1}, if $D = D_1 \by D_2$, then $Z = Z_1 \oplus Z_2$ and $Z_1, Z_2$ are Jordan triple ideals in $Z$. Recall that a subspace $W$ of $Z$ is called a {\em Jordan triple ideal} if $\set{x y z}, \set{y x z} \in W$, for all $x \in W, \ y, z \in Z$. The domain $D$ is of tube--type if and only if each of $D_1, D_2$ is of tube--type. Any bounded symmetric domain $D \subseteq Z$ is the direct product of a finite number of irreducible bounded symmetric domains and $Z$ is the direct sum of their related PJTs. We will be interested only in tube--type domains. We give a list of all irreducible tube--type domains. Any tube--type domain is the direct product of a finite number of the following domains. In the following $n$ will denote a positive integer and $z'$ is the transpose of the square matrix $z \in \BM_n(\IC)$.

\begin{description}
  \item[I$_n$:] $Z = \BM_n(\IC)$, $D = \set{z \in Z \mid \norm{z} < 1}$.
  \item[II$_{2n}$:] $Z = \set{z \in \BM_{2n}(\IC) \mid z' = -z}$, $D = \set{z \in Z \mid \norm{z} < 1}$.
  \item[III$_n$:] $Z = \set{z \in \BM_n(\IC) \mid z' = z}$, $D = \set{z \in Z \mid \norm{z} < 1}$.
  \item[IV$_n$:] $Z = \IC^n$, $n \geq 3$, $D = \set{z \in Z \mid z \cdot \bar{z} + \sqrt{(z \cdot \bar{z})^2 - |z\cdot z|^2} < 2}$. Here $\bar{z}$ is complex conjugate of $z$ and $z \cdot w = \sum_{j = 1}^n z_j w_j$.
  \item[VI:] The exceptional tube--type domain in $Z = \IC^{27}$.
\end{description}
The Jordan triple product related to the PJT structure on $Z$ in the first three cases is $\set{x y z} = (x y^* z + z y^* x)/2$, where $y^*$ is the adjoint of the matrix $y$. In the fourth case, $\set{x y z} = (x \cdot \bar y) z - (x \cdot z) \bar y + (z \cdot \bar y) x$.

The bounded symmetric domains $\mathbf{I}_{n,m}, n \neq m$, $\mathbf{II}_{2n+1}$ and $\mathbf{V}$ are not of tube--type; see \cite[Sec. 4.13]{LOO1}.

A {\em tripotent} is an element $e \in Z$ such that $\set{e e e} = e$. We say that the  two tripotents $e_1, e_2 \in Z$ are {\em orthogonal} if $\set{e_1 e_1 e_2} = 0$. A non--zero tripotent is called {\em primitive} if it is not the sum of non--zero orthogonal tripotents. If $e \in Z$ is a primitive tripotent and $Z$ is the direct sum of the Jordan triple ideals $Z_1, Z_2$, then either $e \in Z_1$ or $e \in Z_2$ and the primitive tripotents in $Z_1$ and $Z_2$ are primitive in $Z$ too. Any non--zero tripotent is the sum of a finite number of pairwise orthogonal primitive tripotents and any primitive tripotent is contained in a maximal set of pairwise orthogonal primitive tripotents. A maximal set of pairwise orthogonal primitive tripotents is called a {\em frame}. All frames in $Z$ have the same number of elements and this number is called the {\em rank} of $Z$. The {\em rank} of a bounded symmetric domain is the rank of the related PJT. A {\em maximal} tripotent is the sum of the elements of some frame. The rank in the cases $\mathbf{I}_{n}$, $\mathbf{II}_{2n}$ and $\mathbf{III}_{n}$ is $n$, the tripotents are exactly the partial isometries and the maximal tripotents are exactly the unitaries. In $\mathbf{IV}_{n}$ the rank is 2 and in $\mathbf{VI}$ the rank is 3.

A {\em generic norm} on $Z$ is a homogeneous polynomial $N: Z \to \IC$ with degree equal to the rank of $Z$ such that $z \in Z$ is invertible if and only if $N(z) \neq 0$ and $N(e) = 1$, for some tripotent $e \in Z$. Thus, if $r$ is the rank of $Z$, then $N(\la z) = \la^r N(z)$, for all $\la \in \IC, \ z \in Z$. The PJT related to an irreducible bounded symmetric domain has a generic norm if and only if the domain is of tube--type. Moreover, any two generic norms differ by a complex factor of absolute value 1. In the case of $\mathbf{I}_{n}$ and $\mathbf{III}_{n}$ the determinant polynomial is a generic norm, and in $\mathbf{II}_{2n}$ the Pfaffian polynomial is a generic norm. In $\mathbf{IV}_{n}$, the polynomial $z \mapsto z \cdot z$ is a generic norm.

Let $D$ be a tube--type domain in $Z$. We can view the polynomial algebra $\CP(Z)$ as a subset of $C(\overline D)$. The Shilov boundary $S$ of $D$ is the Shilov boundary of $\overline D$ relative to the function algebra $\overline{\CP(Z)}$. By \cite[Theorem 6.5]{LOO1}, $S$ is the set of all maximal tripotents of $Z$. Thus, if $\la \in \IT$ and $e \in S$, then $\la e \in S$.

The following proposition describes the structure of the Shilov boundary of a reducible domain. In fact, this elementary result is valid for all bounded symmetric domains. We couldn't find it in the literature, hence, we state it and we prove it here.

\begin{prop} \label{lbl5012}
If a tube--type domain is the direct product of irreducible domains, then its Shilov boundary is the direct product of the Shilov boundaries of the irreducible domains.
\end{prop}

\begin{prf}
It is sufficient to prove that if the tube--type domain $D$ is the direct product of the domains $D_1, D_2$, the Shilov boundary of $D$ equals the product of the Shilov boundaries of $D_1, D_2$. Let $Z, Z_1, Z_2$ be the spaces with Jordan triple systems related to $D, D_1, D_2$, respectively. Then $Z = Z_1 \oplus Z_2$. Let $S, S_1, S_2$ be the Shilov boundaries of $D, D_1, D_2$, respectively and let $r, r_1, r_2$ be the ranks of the domains $D, D_1, D_2$, respectively. Then $r = r_1 + r_2$. If $e \in S$, then $e$ is a maximal tripotent of $Z$. Thus, $e = e_1 + \cdots + e_r$, where $e_1, \cdots, e_r$ are primitive orthogonal tripotents of $Z$. Each of $e_1, \cdots, e_r$ belongs to either $Z_1$ or $Z_2$. Hence, exactly $r_1$ elements of $e_1, \cdots, e_r$ belong to $Z_1$ and their sum is a maximal tripotent in $Z_1$. The other $r_2$ elements of $e_1, \cdots, e_r$ are in $Z_2$ and their sum is a maximal tripotent in $Z_2$. Thus, $e \in S_1 \by S_2$. Conversely, if $e \in S_1, d \in S_2$, then $e = e_1 + \cdots + e_{r_1}$ and $d = d_1 + \cdots + d_{r_2}$, where $e_1, \cdots, e_{r_1}$ are primitive orthogonal tripotents of $Z_1$ and $d_1, \cdots, d_{r_2}$ are primitive orthogonal tripotents of $Z_2$. Thus, $e_1, \cdots, e_{r_1}, d_1, \cdots, d_{r_2}$ is a frame in $Z$ and $e_1 + \cdots + e_{r_1} + d_1 + \cdots + d_{r_2}$ is a maximal tripotent in $Z$. Hence, $e + d$ in $S$.
\end{prf}

Let $D$ be an irreducible tube--type domain with Shilov boundary $S$ and let $N$ be a fixed generic norm on $Z$. Then $|N(z)| = 1$ for all $z \in S$. The {\em reduced Shilov boundary} of the irreducible domain $D$ related to $N$ is defined as the set $S_N = \set{z \in S \mid N(z) = 1}$. In the following proposition we give the properties of the reduced Shilov boundary we need in proving the main theorem in the following section.

\begin{prop} \label{lbl5012a}
Let $D \subseteq Z$ be an irreducible tube--type domain, let $N$ be a generic norm on $Z$ and let $S_N$ be the related reduced Shilov boundary.
\begin{enumerate}
\item $S_N$ is linearly isomorphic to the reduced Shilov boundary of $D$ related to any other generic norm on $Z$.

\item $S_N$ is compact, connected and simply connected.

\item $S_N$ is path connected and locally path connected.
\end{enumerate}
\end{prop}

\begin{prf}
$(1)$ Suppose that $N_1$ is a generic norm on $Z$ and $S_{N_1}$ is the related reduced Shilov boundary. Let $r$ be the rank of $Z$ and let $\la \in \IT$ such that $N = \la N_1$. Then the linear isomorphism $z \mapsto \kappa z$ maps $S_N$ onto $S_{N_1}$, where $\la = \kappa^{r}$.

$(2)$ It is obvious that $S_N$ is compact. By \cite[p. 63]{UPM1}, $S_N$ is a connected and simply connected manifold.

$(3)$ Follows from the fact that $S_N$ is a connected manifold.
\end{prf}

%

\section[Non--Vanishing Functions]{Non--vanishing Functions on the Shilov \\ Boundary of a Tube--Type Domain}
\label{sec2}

Let $n$ be a fixed integer and let $D = D_1 \by \cdots \by D_n$, where $D_1, \dots, D_n$ are irreducible tube--type domains. Suppose that $S_1, \dots, S_n$ are the Shilov boundaries of the domains $D_1, \dots, D_n$, respectively. By Proposition \ref{lbl5012}, $S = S_1 \by \cdots \by S_n$ is the Shilov boundary of $D$. For each $1 \leq j \leq n$, let $Z_j$ be the PJT associated to $D_j$ and let $r_j$ be its rank. Let $N_j$ be a fixed generic norm on $Z_j$ and let $S_{N_j}$ be the related reduced Shilov boundary.

Define the product space $Y = \prod_{j = 1}^n ( \IR \by S_{N_j} )$. To simplify notation, the point $\big( (t_1, u_1), \dots, (t_n, u_n) \big )$ in $Y$ will be denoted by $(t_j, u_j)_{j = 1}^n$. The mapping
\[
\Psi : Y \to S, \ \ (t_j, u_j)_{j = 1}^n \mapsto (e^{\im t_j} u_j)_{j = 1}^n = (e^{\im t_1} u_1, \dots, e^{\im t_n} u_n)
\]
is continuous and surjective. If $\Psi \big( (t_j, u_j)_{j = 1}^n \big) = \Psi \big( (t'_j, u'_j)_{j = 1}^n \big)$, then we write $(t_j, u_j)_{j = 1}^n \sim (t'_j, u'_j)_{j = 1}^n$. This is equivalent to the condition $t'_j = t_j + m_j \al_j, \ u'_j = e^{-m_j \al_j} u_j$, where $m_j \in \IZ$ and $\al_j = 2 \pi / r_j$ for $1 \leq j \leq n$. Note that $\sim$ is an equivalence relation on $Y$. Let $Q$ be the quotient space of $Y$ by this equivalence relation. Note that $Q$ is a compact space since it is the image of the compact subspace $\prod_{j = 1}^n \big( [0, \al_j] \by S_{N_j} \big) \subseteq Y$ under the quotient mapping $q: Y \to Q$.

Define the mapping $\hat \Psi: Q \to S$ by
\[
\big[ (t_j, u_j)_{j = 1}^n \big] \mapsto \Psi \big( (t_j, u_j)_{j = 1}^n \big) = ( e^{\im t_j} u_j)_{j = 1}^n,
\]
where $\big[ (t_j, u_j)_{j = 1}^n \big]$ is the equivalence class of $(t_j, u_j)_{j = 1}^n$. The mapping $\hat \Psi$ is a well--defined bijection and $\Psi = \hat \Psi \circ q$. By the definition of the quotient topology, $\hat \Psi$ is continuous. Therefore, $\hat \Psi$ is a homeomorphism.

In addition to the properties of the space $Y$ given above, we need the following lemma.

\begin{lem} \label{lbl5014}
Let $\nu: Y \to \IT$ be a continuous mapping. Then there exists a continuous mapping $\psi: Y \to \IR$ such that $\nu = e^{\im \psi}$.
\end{lem}

\begin{prf}
By Proposition \ref{lbl5012a}, the space $Y$ is simply connected, path connected and locally path connected. Hence, the fundamental group of $Y$ is the trivial group. Consider the covering space $\IR$ over the space $\IT$ relative to the mapping $p: \IR \to \IT, \ \ t \mapsto e^{\im t}$. By applying the general lifting lemma \cite[Chap. 8, Lemma 14.2]{MUN1}, there exists a continuous function $\psi: Y \to \IR$ such that $\nu = p \circ \psi$. Hence, $\nu = e^{\im \psi}$.
\end{prf}

If $\psi_0, \psi_1: Y \to \IR$ are continuous functions such that $e^{\im \psi_0} = e^{\im \psi_1}$, then there is a unique integer $k$ such that $\psi_0 = \psi_1 + 2 \pi k$. This remark will be used frequently in the proof of the following lemma.

\begin{lem} \label{lbl5016}
If $\vfi: S \to \IT$ is continuous, then $\vfi = \theta e^{\im \psi}$ for some continuous function $\psi: S \to \IR$, where $\theta: S \to \IC$ is the unimodular function defined by the unique integers $k_1, \dots, k_n$ as
\[
\theta \big( (u_j)_{j = 1}^n \big) = \prod_{j = 1}^n N_j^{k_j}(u_j),  \quad (u_j)_{j = 1}^n \in S.
\]
\end{lem}

\begin{prf}
We denote the Kronecker delta function by $\delta_{jl}$. Consider the continuous mapping $\nu_0: Y \to \IT$, where $\nu_0 = \vfi \circ \Psi$. By Lemma \ref{lbl5014}, there exists a continuous function $\psi_0: Y \to \IR$ such that $\nu_0 = e^{\im \psi_0}$. If $l$ is an integer such that $1 \leq l \leq n$, then
{\setlength\arraycolsep{2pt}
\begin{eqnarray*}
  \exp \Big( \im \psi_0 \big( (t_j + \delta_{jl} \al_j, e^{-\im \delta_{jl} \al_j}u_j)_{j = 1}^n \big) \Big)
  &=& \nu_0 \big( (t_j + \delta_{jl} \al_j, e^{-\im \delta_{jl} \al_j}u_j)_{j = 1}^n \big) \\
  &=& \vfi \big( (e^{\im t_j} u_j)_{j = 1}^n  \big) \\
  &=& \nu_0 \big ( (t_j, u_j)_{j = 1}^n \big ) \\
  &=& \exp \Big( \im \psi_0 \big( (t_j, u_j)_{j = 1}^n \big) \Big),
\end{eqnarray*}}%
for all $(t_j, u_j)_{j = 1}^n \in Y$. Hence, there exists an integer $k_l$ such that
{\setlength\arraycolsep{2pt}
\begin{eqnarray}
  \psi_0 \big( (t_j + \delta_{jl} \al_j, e^{-\im \delta_{jl} \al_j}u_j)_{j = 1}^n \big) = \psi_0 \big( (t_j, u_j)_{j = 1}^n \big) + 2 \pi k_l.  \label{eqn5002} %
\end{eqnarray}}%

Using the $n$ integers $k_1, \dots, k_n$ satisfying Eq. (\ref{eqn5002}), we define the unimodular function $\theta:S \to \IC$ such that $\theta \big( (u_j)_{j = 1}^n \big) = \prod_{j = 1}^n N_j^{k_j}(u_j)$ for $(u_j)_{j = 1}^n \in S$.

Define the function $\vfi_1: S \to \IT$ such that $\vfi_1 = \theta^{-1} \vfi$ and let ${\nu_1: Y \to \IT}$ be the function $\nu_1 = \vfi_1 \circ \Psi$. By Lemma \ref{lbl5014}, there exists a continuous function $\psi_1: Y \to \IR$ such that $\nu_1 = e^{\im \psi_1}$.

If $(t_j, u_j)_{j = 1}^n \in Y$, then
{\setlength\arraycolsep{2pt}
\begin{eqnarray*}
\exp \Big( \im \psi_1 \big( (t_j, u_j)_{j = 1}^n \big) \Big) &=& \nu_1 \big( (t_j, u_j)_{j = 1}^n \big) = \vfi_1 \big( (e^{\im t_j} u_j)_{j = 1}^n \big) \\
&=& \Big(\theta\big( (e^{\im t_j} u_j)_{j = 1}^n \big)\Big)^{-1} \vfi \big( (e^{\im t_j} u_j)_{j = 1}^n \big) \\
&=& N_1^{-k_1}(e^{\im t_1} u_1) \cdots N_n^{-k_n}(e^{\im t_n} u_n) \vfi \big( (e^{\im t_j} u_j)_{j = 1}^n \big) \\
&=& e^{- \im \sum_{j = 1}^n r_j k_j t_j} \nu_0 \big( (t_j, u_j)_{j = 1}^n \big) \\
&=& \exp \left( \im \psi_0 \big( (t_j, u_j)_{j = 1}^n \big) - \im \sum_{j = 1}^n r_j k_j t_j \right).
\end{eqnarray*}}%
Hence, there is an integer $k$ such that
{\setlength\arraycolsep{2pt}
\begin{eqnarray}
\psi_1 \big( (t_j, u_j)_{j = 1}^n \big) = \psi_0 \big( (t_j, u_j)_{j = 1}^n \big) - \sum_{j = 1}^n r_j k_j t_j + 2 \pi k, \label{eqn5003} %
\end{eqnarray}}%
for all $(t_j, u_j)_{j = 1}^n \in Y$. Let $l$ be an integer such that $1 \leq l \leq n$. Then, by Eq. (\ref{eqn5003}), we get
{\setlength\arraycolsep{2pt}
\begin{eqnarray*}
\psi_1 \big( (t_j + \delta_{jl} \al_j, e^{-\im \delta_{jl} \al_j}u_j)_{j = 1}^n \big) &=&
   \psi_0 \big( (t_j + \delta_{jl} \al_j, e^{-\im \delta_{jl} \al_j}u_j)_{j = 1}^n \big)\\
 & & - \sum_{j = 1}^n r_j k_j t_j - r_l k_l \al_l + 2 \pi k.%
\end{eqnarray*}}%
Now we use Eq. (\ref{eqn5002}) and the equality $\al_l = 2 \pi / r_l$ to get
{\setlength\arraycolsep{2pt}
\begin{eqnarray*}
 \psi_1 \big( (t_j + \delta_{jl} \al_j, e^{-\im \delta_{jl} \al_j}u_j)_{j = 1}^n \big) &=& \psi_0 \big( (t_j, u_j)_{j = 1}^n \big) + 2 \pi k_l \\
 & & - \sum_{j = 1}^n r_j k_j t_j - r_l k_l \al_l + 2 \pi k \\
 &=& \psi_0 \big( (t_j, u_j)_{j = 1}^n \big) - \sum_{j = 1}^n r_j k_j t_j + 2 \pi k \\
 &=& \psi_1 \big( (t_j, u_j)_{j = 1}^n \big).   %
\end{eqnarray*}}%
Note that we obtained the last equality from Eq. (\ref{eqn5003}). By induction, it follows that the relation
{\setlength\arraycolsep{2pt}
\begin{eqnarray}
\psi_1 \big( (t_j + m_j \al_j, e^{-\im m_j \al_j} u_j)_{j = 1}^n \big) = \psi_1 \big( (t_j, u_j)_{j = 1}^n \big) \label{eqn5006} %
\end{eqnarray}}%
holds for all $(t_j, u_j)_{j = 1}^n \in Y$ and all $m_1, \dots, m_n \in \IZ$.

Let $Q$ be the quotient space defined above. Define the function
\[
\hat \psi_1: Q \to \IR, \quad \big[ (t_j, u_j)_{j = 1}^n \big] \mapsto \psi_1 \big( (t_j, u_j)_{j = 1}^n \big).
\]
By Eq. (\ref{eqn5006}) and by the remarks on the definition of the equivalence relation $\sim$ on $Y$, the function $\hat \psi_1$ is well--defined and $\psi_1 = \hat \psi_1 \circ q$, where $q$ is the quotient canonical mapping. Moreover, it follows that $\hat \psi_1$ is continuous.

Let $\hat \Psi: Q \to S$ be the homeomorphism defined above and define the continuous function $\psi: S \to \IR$, where $\psi = \hat \psi_1 \circ \hat \Psi^{-1}$. Let $(u_j)_{j = 1}^n \in S$. For all $j = 1, \dots, n$, let $s_j$ be a real number such that $0 \leq s_j < \al_j$ and $N_j(u_j) = e^{\im r_j s_j}$. Hence,
{\setlength\arraycolsep{2pt}
\begin{eqnarray*}
\exp \Big( \im \psi \big( (u_j)_{j = 1}^n \big) \Big) &=& \exp \Big( \im \hat \psi_1 \circ \hat \Psi^{-1} \big ( (u_j)_{j = 1}^n \big ) \Big) \\
&=& \exp \Big( \im \psi_1 \big( (s_j, e^{-\im s_j} u_j)_{j = 1}^n \big) \Big) \\
&=& \nu_1 \big( (s_j, e^{-\im s_j} u_j)_{j = 1}^n \big) \\
&=& \vfi_1 \Big( \Psi \big( (s_j, e^{-\im s_j} u_j)_{j = 1}^n \big) \Big) 
= \vfi_1 \big( (u_j)_{j = 1}^n \big). %
\end{eqnarray*}}%

Thus, $\vfi_1 = e^{\im \psi}$. Since $\vfi_1 = \theta^{-1} \vfi$, it follows that $\vfi = \theta e^{\im \psi}$. This proves the existence part of this lemma.

To prove the uniqueness part, it is sufficient to show that if $k'_1, \dots, k'_n$ are integers and $\tilde \psi : S \to \IR$ is a continuous function, then the equality
\[
\exp \Big(\im \tilde \psi \big( (u_j)_{j = 1}^n \big) \Big) \prod_{j = 1}^n N_j^{k'_j}(u_j) = 1 \text{ for all } (u_j)_{j = 1}^n \in S
\]
implies $k'_1 = \cdots = k'_n = 0$.

If $(t_j, u_j)_{j = 1}^n \in Y$, then
\[
1 = \exp \Big( \im \tilde \psi \big( (e^{\im t_j} u_j)_{j = 1}^n \big) \Big) \prod_{j = 1}^n N_j^{k'_j}(e^{\im t_j} u_j) = \exp \Big( \im \tilde \psi \big( (e^{\im t_j} u_j)_{j = 1}^n \big) + \im \sum\limits_{j=1}^n k'_j r_j t_j \Big).
\]
Hence, there exists an integer $k$ such that
\[
\tilde \psi \big( (e^{\im t_j} u_j)_{j = 1}^n \big) + \sum_{j=1}^n k'_j r_j t_j = 2 \pi k.
\]
Choose a fixed element $u_j \in S_{N_j}$ for each $1 \leq j \leq n$. Let $m_1, \dots, m_n \in \set{0, 1}$. By the previous equation we get
{\setlength\arraycolsep{2pt}
\begin{eqnarray*}
\tilde \psi \big( (u_j)_{j = 1}^n \big) &=& \tilde \psi \big( (e^{\im \cdot 0} u_j)_{j = 1}^n \big) + \sum_{j=1}^n k'_j r_j \cdot 0 = 2 \pi k \\
&=& \tilde \psi \big( (e^{2 \im \pi m_j} u_j)_{j = 1}^n \big) + 2 \pi \sum_{j=1}^n k'_j r_j m_j \\
&=& \tilde \psi \big( (u_j)_{j = 1}^n \big) + 2 \pi \sum_{j=1}^n k'_j r_j m_j.
\end{eqnarray*}}%
Thus, $\sum_{j=1}^n k'_j r_j m_j = 0$. It follows easily that $k'_1 = \cdots = k'_n = 0$.
\end{prf}

The next theorem is the main result in this section, it follows almost directly from Lemma \ref{lbl5016}.

\begin{thm} \label{lbl5018}
For every continuous non--vanishing function $\vfi$ on $S$ there is a continuous function $\psi$ on $S$ such that $\vfi = \theta e^\psi$ where $\theta: S \to \IC$ is the unimodular function defined as $\theta \big( (u_j)_{j = 1}^n \big) = \prod_{j = 1}^n N_j^{k_j}(u_j)$ for the unique integers $k_1, \dots, k_n$.
\end{thm}

\begin{prf}
Let $\vfi_r: S \to \IR$ and $\vfi_u: S \to \IT$ be the continuous functions defined by the relations $\vfi_r = |\vfi|, \vfi_u = \vfi/|\vfi|$, respectively. Since $\vfi_r$ is positive, it is easy to show that there is a continuous function $\psi_r: S \to \IR$ such that $\vfi_r = e^{\psi_r}$. Lemma \ref{lbl5016} asserts that there exist unique integers $k_1, \dots, k_n$ such that $\vfi_u = \theta e^{\im \psi_u}$ where $\theta: S \to \IC$ is the unimodular function defined by $\theta \big( (u_j)_{j = 1}^n \big) = \prod_{j = 1}^n N_j^{k_j}(u_j)$ and $\psi_u: S \to \IR$ is a continuous function. Define the continuous function $\psi : S \to \IC$, where $\psi = \psi_r + \im \psi_u$. Hence, $\vfi = \theta e^\psi$. Since $\vfi_u = \vfi/|\vfi|$, the integers $k_1, \dots, k_n$ are unique.
\end{prf}

\begin{rem} \label{lbl5017a}
The integers $k_1, \dots, k_n$ in Theorem \ref{lbl5018} are independent of the choice of the generic norms $N_1, \dots, N_n$.
\end{rem}

%

\section{Harmonic Analysis on Tube--Type Domains}
\label{sec3}

Throughout this section, we consider a fixed irreducible tube--type domain $D \subseteq Z$ with its related PJT structure and let $S$ be its Shilov boundary. We choose a fixed generic norm $N$ on $Z$. There is a unique $K$--invariant regular Borel probability measure $\mu$ on $S$, where $K$ is the group of all invertible linear mappings on $Z$ such that $g(D) = D$; see \cite[p. 125]{UPM1}.

Let $L^2(S)$ be the space of square--integrable functions on $S$ with respect to~$\mu$. The Hardy space $H^2(S)$ is the $L^2$ closure of the algebra of all polynomials $\CP(Z)$ in $L^2(S)$. We denote the orthogonal projection of the space $L^2(S)$ onto $H^2(S)$ by $P_D$.

Define the set $V_D = \{v \in Z \mid N(v) = 0\}$. Note that the set $V_D$ is independent of the choice of the generic norm $N$. For any $u \in  Z$, we define the homogeneous polynomial $\ell_u$ of degree~$1$ on $Z$ by $\ell_u(z) = \innerp{z}{u}$, where $\innerp{\cdot}{\cdot}$ is the inner product defined in Eq. (\ref{eqn5000}). Let $\CH(Z)$ be the linear subspace of $\CP(Z)$ spanned by the polynomials $\ell_v^k$, where $k$ is a non--negative integer and $v \in V_D$. This space is called the space of {\em harmonic polynomials} \cite[p. 121]{UPM1}. 

Let $S_N$ be the reduced Shilov boundary of $D$ related to $N$.
There exists a unique probability Borel measure on $S_N$ which is invariant under the commutator subgroup of $K$. We define $L^2(S_N)$ as the $L^2$ space relative to this measure. If $N_1$ is another generic norm on $Z$ and $S_{N_1}$ its related reduced Shilov boundary, then, by Proposition \ref{lbl5012a}, there is a linear isomorphism $S_{N_1} \to S_N$. This isomorphism preserves the invariant measures and it induces a unitary $L^2(S_N) \to L^2(S_{N_1})$.

\begin{prop}
There is a collection $\set{H_\al}_{\al \in I}$ of finite--dimensional subspace of $L^2(S)$ for some indexing set $I$ such that:
\begin{enumerate}
\item For each $\al \in I$, the elements of $H_\al$ are harmonic polynomials restricted to $S$.
\item The space $L^2(S_N)$ is unitarily equivalent to the orthogonal sum $\bigoplus_{\al \in I} H_\al$.
\item $L^2(S) = \bigoplus_{\al \in I, \, l \in \IZ} \, N^l H_\al$.
\item $H^2(S) = \bigoplus_{\al \in I, \, l \geq 0} \, N^l H_\al$.
\end{enumerate}
\end{prop}

\begin{prf}
By \cite[Theorem 2.8.68]{UPM1}, the space $L^2(S_N)$ can be identified with the orthogonal sum $\bigoplus_{\al \in I} H_\al$, for some indexing set $I$. For each $\al \in I$, $H_\al$ is a finite--dimensional Hilbert subspace of $L^2(S)$; see \cite[p.118 and p.133]{UPM1}. Moreover, by \cite[Proposition 4.11.51]{UPM1}, the elements of $H_\al$ are harmonic polynomials restricted to $S$. Thus, we have proved $(1)$ and $(2)$. Statement $(3)$ follows from Theorem 2.8.38 and Eq. (2.8.34) in \cite{UPM1}. The equality in $(4)$ is proved in Theorem 2.8.47 in \cite{UPM1}.
\end{prf}

For every $\al \in I$, let the positive integer $d_\al$ be the dimension of $H_\al$ and let the set of harmonic polynomials $\set{p_1^{(\al)}, \ldots, p_{d_\al}^{(\al)} }$ be an orthonormal basis of $H_\al$. Thus, $\set{N^l p_k^{(\al)} \mid l \in \IZ, 1 \leq k \leq d_\al, \al \in I}$ is an orthonormal basis for $L^2(S)$ and $\set{N^l p_k^{(\al)} \mid l \geq 0, 1 \leq k \leq d_\al, \al \in I}$ is an orthonormal basis for $H^2(S)$. Therefore,
\[
P_D(N^l p_k^{(\al)}) = \left\{
\begin{array}{ll}
N^l p_k^{(\al)}, & \textrm{if $l \geq 0,$}\\
0, & \textrm{if $l < 0$.}\\
\end{array} \right.
\]

Now we define the unitary mapping $U: L^2(S) \to L^2(\IT) \ot L^2(S_N)$ such that $U(N^l p_k^{(\al)}) = z^l \ot p_k^{(\al)}$ for $\al \in I, \ 1 \leq k \leq d_\al, \ l \in \IZ$, where $L^2(\IT)$ is the space of square--integrable functions on $\IT$ relative to the normalized arc length measure and $z$ is the inclusion map $z: \IT \to \IC, \ \la \mapsto \la$.

Define the $*$--isomorphism $\pi : \BDD{L^2(S)} \to \BDD{L^2(\IT) \ot L^2(S_N)}$ by setting $\pi(a) = UaU^*$ for $a \in \BDD{L^2(S)}$.

\begin{prop} \label{lbl6001}
If $u \in Z$, then there exists a unique pair of operators $A_u, B_u$ in $\BDD{L^2(S_N)}$ such that $\overline \ell_u p = N^{-1} A_u(p) + B_u(p)$, for all $p \in \CH(Z)$.
\end{prop}

\begin{prf}
The result follows from Eqs. (4.11.50) and (4.11.58) in \cite{UPM1}.
\end{prf}

The mapping which takes a continuous function on $S$ to its related multiplication operator on $L^2(S)$ is a $*$--isomorphism. Thus, we can consider $C(S)$ as a unital $\cstr$--subalgebra of $\BDD{L^2(S)}$. Let $H^2(\IT)$ be the Hardy space related to $L^2(\IT)$ and let $P_\IT$ be the orthogonal projection of $L^2(\IT)$ onto the subspace $H^2(\IT)$.

\begin{prop} \label{lbl6002}
Let $u \in Z$ and let $A_u, B_u$ be the unique pair of operators related to $u$ as in Proposition \ref{lbl6001}. Then we have
\begin{enumerate}
\item $\pi(N) = z \ot 1$.

\item $\pi(P_D) = P_\IT \ot 1$.

\item $\pi(\ell_u) = z \ot A_u^* + 1 \ot B_u^*$.  
\end{enumerate}
\end{prop}

\begin{prf}
Let $z^l \ot p_k^{(\al)}$ be an arbitrary element in the orthonormal basis of $L^2(\IT) \ot L^2(S_N)$ defined above. Then
{\setlength\arraycolsep{2pt}
\begin{eqnarray*}
  \pi(N) (z^l \ot p_k^{(\al)}) &=& U N U^* (z^l \ot p_k^{(\al)})
                                   = U \big( N  (N^l p_k^{(\al)}) \big) \\
                                  &=& U (N^{l + 1} p_k^{(\al)})
                                   = z^{l + 1} \ot p_k^{(\al)} \\
                                  &=& (z \ot 1) (z^l \ot p_k^{(\al)}).
\end{eqnarray*}}%
Thus, $\pi(N) = z \ot 1$. This proves $(1)$.

To show that $(2)$ holds, note first that $\pi(P_D) (z^l \ot p_k^{(\al)}) = U \big( P_D (N^l p_k^{(\al)}) \big)$. Thus, if $l \geq 0$, then
{\setlength\arraycolsep{2pt}
\begin{eqnarray*}
\pi(P_D) (z^l \ot p_k^{(\al)}) &=& U(N^l p_k^{(\al)}) = z^l \ot p_k^{(\al)} \\
                               &=& P_\IT(z^l) \ot p_k^{(\al)} = (P_\IT \ot 1) (z^l \ot p_k^{(\al)}).
\end{eqnarray*}}%
Similarly, if $l < 0$, then $\pi(P_D) (z^l \ot p_k^{(\al)}) = 0 = (P_\IT \ot 1) (z^l \ot p_k^{(\al)})$. This proves~$(2)$.

Now we prove $(3)$. By Proposition \ref{lbl6001} above,
{\setlength\arraycolsep{2pt}
\begin{eqnarray*}
  \pi(\overline \ell_u) (z^l \ot p_k^{(\al)}) &=& U \big( \overline \ell_u (N^l p_k^{(\al)}) \big) = U (N^l \overline \ell_u p_k^{(\al)}) \\
   &=& U \Big( N^l \big( N^{-1} A_u(p_k^{(\al)}) + B_u(p_k^{(\al)}) \big) \Big) \\
   &=& U \big( N^{l-1} A_u (p_k^{(\al)}) + N^l B_u (p_k^{(\al)}) \big) \\
   &=& z^{l-1} \ot A_u (p_k^{(\al)}) + z^l \ot B_u (p_k^{(\al)}) \\
   &=& (z^{-1} \ot A_u + 1 \ot B_u)(z^l \ot p_k^{(\al)}).
\end{eqnarray*}}%
Thus, $\pi(\overline \ell_u) = z^{-1} \ot A_u + 1 \ot B_u$. Therefore, $\pi(\ell_u) = z \ot A_u^* + 1 \ot B_u^*$.
\end{prf}

Now we define $\CC_D$ as the unital $\cstr$--subalgebra of $\BDD{L^2(S_N)}$ generated by the operators $A_u, B_v$ for all $u, v \in Z$. We need the following property of the $\cstr$--algebra $\CC_D$ which follows easily from \cite[Proposition 4.11.140]{UPM1}.

\begin{prop} \label{lbl6003}
There exists a character $\tau$ on the unital $\cstr$--algebra $\CC_D$.
\end{prop}

\begin{prop} \label{lbl6004}
The $*$--isomorphism $\pi$ maps $C(S)$ into $C(\IT) \ot \CC_D$.
\end{prop}

\begin{prf}
Let $C_0$ be the unital $*$--subalgebra of $C(S)$ generated by all of the functions $\ell_u$, where $u \in Z$. Let $u_1, u_2 \in S$ such that $u_1 \neq u_2$. If we set $v = u_1 - u_2$, then $\ell_v(u_1) \neq \ell_v(u_2)$. Thus, $C_0$ separates the points of $S$ and by the Stone--Weierstrass theorem, $C_0$ is dense in $C(S)$. By Proposition \ref{lbl6002}, $\pi(\ell_u) = z \ot A_u^* + 1 \ot B_u^* \in C(\IT) \ot \CC_D$, for all $u \in Z$. Therefore, $\pi(C_0) \subseteq C(\IT) \ot \CC_D$ and the result follows by density.
\end{prf}

\section{Toeplitz Operators and Their Index Theorem}
\label{sec4}

We keep the notation of the previous section. If $\vfi \in C(S)$, then we define the Toeplitz operator $W_\vfi$ with symbol $\vfi$ as
\[
W_\vfi(g) = P_D(\vfi g) \quad (g \in H^2(S)).
\]
We define the Toeplitz algebra $\CT_D$ related to the tube-type domain $D$ as the $\cstr$-algebra generated by all these Toeplitz operators. In this section we prove an index theorem for these operators using the general index theory introduced by Murphy in~\cite{MUR1}.

Let $B$ be a unital $\cstr$--algebra. An {\em indicial triple} for $B$ is an ordered triple $\OM = (\CL, F, \tr)$ such that
\begin{itemize}
\item $\CL$ is a unital $\cstr$--algebra and $B$ is a unital $\cstr$--subalgebra of $\CL$,
\item $F$ is a self--adjoint unitary in $\CL$,
\item $\tr$ is a lower semicontinuous trace on $\CL$,
\item the set $\set{\vfi \in B \mid [F, \vfi] = F \vfi - \vfi F \in \CM_\tr}$ is dense in $B$, where $\CM_\tr$ is the definition ideal of the trace $\tr$.
\end{itemize}

We will denote the closure of $\CM_\tr$ by $\CK_\tr$. The commutator $[F, \vfi]$ will always be denoted by $\df \vfi$. The unital dense $*$--subalgebra $\set{\vfi \in B \mid \df \vfi \in \CM_\tr}$ will be denoted by $B_\OM$. Let $\CL_\IT$ be the $\cstr$--subalgebra of $\BDD{L^2(\IT)}$ generated by $C(\IT)$ and the projection $P_\IT$. Then the triple $\OM_\IT = (\CL_\IT, F_\IT, \tr_\IT)$ is an indicial triple for $C(\IT)$, where $F_\IT = 2P_\IT - 1$ and $\tr_\IT$ is the restriction of the canonical trace on $\BDD{L^2(\IT)}$ to $\CL_\IT$; see Example 2.7 in \cite{MUR1}.

By Proposition \ref{lbl6004}, $\pi\big(C(S)\big) \subseteq \CL_\IT \ot \CC_D$. We denote the unital $\cstr$--algebra $\CL_\IT \ot \CC_D$ by $\CL_D$. From now on, we will ignore the $*$--isomorphism $\pi$ and we will consider that $C(S)$ is a unital $\cstr$--subalgebra of $\CL_D$. By Proposition \ref{lbl6002}, $P_D \in \CL_D$.

Define the lower semicontinuous $\tr_\IT \ot \tau$ on $\CL_D$ by the relation
\[
(\tr_\IT \ot \tau)(a) = \tr_\IT \big((\iota \ot \tau)(a) \big) \text{ for } a \in \CL_D^+,
\]
where $\tau$ is the character on the $\cstr$--algebra $\CC_D$ in Proposition \ref{lbl6003}. Let $a \in \CM_{\tr_\IT}$. Then, by \cite[Proposition A1(a)]{PHRA1}, $|a| \in \CM_{\tr_\IT}^+$. Thus, $|a \ot b| = |a| \ot |b| \in \CM_{\tr_\IT \ot \tau}^+$ for all $b \in \CC_D$. Again by the same proposition we get $a \ot b \in  \CM_{\tr_\IT \ot \tau}$. By density it follows that $\CK_{\tr_\IT} \ot \CC_D \subseteq \CK_{\tr_\IT \ot \tau}$. Now we define the lower semicontinuous trace $\tr_D$ on $\CL_D$ by the equality
\[
\tr_D(a) = \left\{ \begin{array}{ll}
(\tr_\IT \ot \tau)(a), & a \in (\CK_{\tr_\IT} \ot \CC_D)^+,\\
+\infty,  & a \in \CL_D^+ \backslash (\CK_{\tr_\IT} \ot \CC_D)^+.\\
\end{array} \right.
\]
Moreover, $\CK_{\tr_D} = \CK_{\tr_\IT} \ot \CC_D$. The following proposition shows that the trace $\tr_D$ defines an indicial triple for $C(S)$.

\begin{prop} \label{lbl6008}
Let $F_D = 2P_D - 1$. Then the triple $\OM_D = (\CL_D, F_D, \tr_D)$ is an indicial triple for $C(S)$.
\end{prop}

\begin{prf}
Let $B_{\OM_D}$ be the unital $*$--subalgebra of all elements $\vfi \in C(S)$ such that $\df \vfi \in \CM_{\tr_D}$. By Example 2.7 in \cite{MUR1}, $\df \ell_u = [F_\IT, z] \ot A_u^* \in \CM_{\tr_\IT} \ot \CC_D$ for all $u \in Z$.  Thus, $\df \ell_u \in \CM_{\tr_D}$. This implies that $\ell_u \in B_{\OM_D}$ for all $u \in Z$. Hence, $B_{\OM_D}$ contains the unital dense $*$--subalgebra $C_0$ defined in the proof of Proposition \ref{lbl6004}. This implies that $B_{\OM_D}$ is dense in $C(S)$. Therefore, $\OM_D$ is an indicial triple for $C(S)$.
\end{prf}

A {\em topological index} on a unital $\cstr$--algebra $B$ is a locally constant homomorphism from $\Inv(B)$, the group of invertible elements in $B$, to the additive group $\IR$. By \cite[Theorem 2.6]{MUR1}, if $\OM = (\CL, F, \tr)$ is an indicial triple for $B$, then there is a unique topological index $\om$ on $B$ such that $\om(\vfi) = \frac12 \tr(\vfi^{-1} \df \vfi)$, for all $\vfi \in \Inv(B_\OM)$.

Let $\om_D$ be the unique topological index on $C(S)$ related to the indicial triple $\OM_D$.

\begin{prop} \label{lbl6008b}
Let $\vfi$ be a non-vanishing function on $S$. Then $\om_D(\vfi) = k$, where $k$ is the unique integer such that $\vfi = N^k e^\psi$ for some $\psi \in C(S)$.
\end{prop}

\begin{prf}
First we calculate $\om_D(N)$,
{\setlength\arraycolsep{2pt}
\begin{eqnarray*}
\om_D(N) &=& \frac12 \tr_D(N^{-1} \df N) = \frac12 \tr_D \big( (z^{-1} \ot 1) (\df z \ot 1) \big) \\
       &=& \frac12 (\tr_\IT \ot \tau) \big( (z^{-1} \df z) \ot 1 \big) = \frac12 \tr_\IT(z^{-1} \df z) \tau(1) = \om_\IT(z),
\end{eqnarray*}}%
where $\om_\IT$ is the unique topological index on $C(\IT)$ related to the indicial triple $\OM_\IT$. Hence, by Example 2.7 in \cite{MUR1}, $\om_D(N) = 1$. By Theorem \ref{lbl5018}, there exists a unique integer $k$ such that $\vfi = N^k e^\psi$, where $\psi \in C(S)$. By elementary properties of topological indices, we have
\[
\om_D(\vfi) = \om_D(N^k) + \om_D(e^\psi) = k \cdot 1 + 0 = k.
\]
Thus, $\om_D(\vfi) = k$.
\end{prf}

It is obvious that the value of the topological index $\om_D$ is independent of the choice of the character $\tau$. Moreover, by Remark \ref{lbl5017a}, it is independent of the choice of the generic norm $N$.

Suppose that $\OM = (\CL, F, \tr)$ is an indicial triple for $B$. Let $\vfi \in B$. The element $T_\vfi = P \vfi P$ is called a {\em Toeplitz element} associated to $\OM$, where $P = (F + 1)/2$. Let $\CA$ be the $\cstr$--subalgebra generated by the Toeplitz elements associated to $\OM$. This unital algebra is called the {\em Toeplitz algebra} associated to $\OM$. An element $a \in \CA$ is called an {\em $\OM$--Fredholm element} if it is invertible modulo $\CK_\Tr$, where the trace $\Tr$ is the restriction of $\tr$ to $\CA$. If $a \in \CA$ is an $\OM$--Fredholm element, then there exists $b \in \CA$ such that $P - ab, P - ba \in \CM_\Tr$ and the {\em $\OM$--Fredholm index} of $a$ is defined as $\ind(a) = \Tr(ab - ba)$. The $\OM$--Fredholm index is a well--defined function. This function is real--valued and it has algebraic properties similar to the classical Fredholm index; see \cite{MUR4, MUR1}. Now by Proposition \ref{lbl6008b} and \cite[Theorem 3.1]{MUR1}, we have the following immediate result.

\begin{thm} \label{lbl6009}
Let $\vfi$ be a non--vanishing function on $S$. Then $T_\vfi$ is an $\OM_D$--Fredholm element and $\ind_D(T_\vfi) = -k$, where $k$ is the unique integer such that $\vfi = N^k e^\psi$ for some $\psi \in C(S)$ and $\ind_D$ is the $\OM_D$--Fredholm index.
\end{thm}

For any $\cstr$--algebra $A$, we denote the closed ideal generated by all of the commutators in $A$ by $\Com(A)$ and we will call it the {\em closed commutator ideal} of $A$. We denote the ideal of compact operators on $H^2(\IT)$ by $\CPT{H^2(\IT)}$.

Let $\CA_D$ be the Toeplitz algebra related to the indicial triple $\OM_D$. The algebra $\CA_D$ is $*$--isomorphic to the Toeplitz algebra $\CT_D$ of Toeplitz operators on the tube--type domain $D$ defined in the beginning of this section and the $*$--isomorphism $\CT_D \to \CA_D$ is defined by $a \mapsto aP_D$. By \cite[Theorem 4.11.76]{UPM1}, $\Com(\CT_D) = \CPT{H^2(\IT)} \ot \CC_D$. The latter $*$--isomorphism maps this ideal onto $\Com(\CA_D) = P_\IT \CPT{L^2(\IT)} P_\IT \ot \CC_D$. By Lemma 3.1 in \cite{ADBA1}, $\CPT{L^2(\IT)} = \CK_{\tr_\IT}$. Hence, $\Com(\CA_D) = P_\IT \CK_{\tr_\IT} P_\IT \ot \CC_D$.

\begin{thm} \label{lbl6011}
If $\vfi \in C(S)$, then the Toeplitz element $T_\vfi$ is $\OM_D$--Fredholm if and only if $\vfi$ is invertible.
\end{thm}

\begin{prf}
If $a \in \CK_{\Tr_D}$, then $a \in \CK_{\tr_D} = \CK_{\tr_\IT} \ot \CC_D$. By density, we can choose a sequence $(a_n)_{n \geq 1}$ in the algebraic tensor product $\CK_{\tr_\IT} \odot \CC_D$ such that $a_n \to a$. Thus, $P_D a_n P_D \in P_\IT \CK_{\tr_\IT} P_\IT \odot \CC_D$. Therefore, $a = P_D a P_D \in P_\IT \CK_{\tr_\IT} P_\IT \ot \CC_D$, since $P_D a_n P_D \to P_D a P_D$. Thus, $\CK_{\Tr_D} \subseteq \Com(\CA_D)$.

By \cite[Theorem 3.11]{UPM3}, there exists a $*$--homomorphism $\ro : \CA_D \to C(S)$ such that $\ro(T_\vfi) = \vfi$, for all $\vfi \in C(S)$, and $\ker(\ro) = \Com(\CA_D)$. Now if $T_\vfi$ is an $\OM_D$--Fredholm element, then there is an element $S \in \CA_D$ such that both of the elements $P_D - S T_\vfi$ and $P_D - T_\vfi S$  belong to the ideal $\CK_{\Tr_D}$ and it follows immediately that $\ro(S) \vfi = \vfi \ro(S) = 1$. Hence, $\vfi$ is an invertible function. The opposite direction follows from Theorem \ref{lbl6009}.
\end{prf}

\begin{cor} \label{lbl6012}
If $\vfi$ be a non--vanishing function on $S$, then $\ind_\OM (T_\vfi) = 0$ if and only if $\vfi = e^\psi$, for some $\psi \in C(S)$.
\end{cor}

Now we extend our results to Toeplitz operators with $M \by M$ matrix symbols defined on the Shilov boundary of an irreducible tube--type domain. These operators act on the Hilbert space $H^2_M(S) = \underbrace{H^2(S) \oplus \cdots \oplus H^2(S)}_{M \text{ times}}$. The space $L^2_M(S)$ is defined similarly and we denote the orthogonal projection of $L^2_M(S)$ onto $H^2_M(S)$ by $P_M$. Let $C_M$ be the space of continuous $M \by M$ matrix--valued functions on $S$. Then $C_M = \MTRX_M\big( C(S) \big)$. The Toeplitz operator with symbol $\Phi \in C_M$ is defined by $T_\Phi (\mathbf{f}) = P_M(\Phi \mathbf{f})$ for all $\mathbf{f} \in H^2_M(S)$. These Toeplitz operators can be obtained from the indicial triple $\OM_M = (\CL_M, F_M, \tr_M)$ for $C_M$, where $\CL_M = \MTRX_M(\CL_D) = \CL_D \ot \MTRX_M(\IC)$, $F_M = F_D \ot 1_M = 2P_M - 1$ and the trace $\tr_M$ is defined by $\tr_M([a_{jk}]) = \sum_{j = 1}^M \tr_D(a_{jj})$, $[a_{jk}] \in \CL_M^+$; see \cite[Theorem 2.9] {MUR1}. The following is the index theorem for Toeplitz operators with $M \by M$ matrix symbols on an irreducible tube--type domain. The proof of the following theorem is similar to that of Theorem 3.7 in \cite{ADBA1}, hence, we will not write it.

\begin{thm} \label{lbl6012a}
Let $\Phi \in C_M$. Then $T_\Phi$ is an $\OM_M$--Fredholm operator if and only if $\det \Phi$ is non-vanishing on $S$ and in this case we have
\[
\ind_M(T_\Phi) = -k,
\]
where $k$ is the unique integer such that $\det \Phi = N^k e^\psi$, for some continuous function $\psi$ on $S$ and $\ind_M$ is the $\OM_M$--index mapping.
\end{thm}

%

\bibliographystyle{amsplain}
\providecommand{\bysame}{\leavevmode\hbox to3em{\hrulefill}\thinspace}
\providecommand{\MR}{\relax\ifhmode\unskip\space\fi MR }
\providecommand{\MRhref}[2]{%
  \href{http://www.ams.org/mathscinet-getitem?mr=#1}{#2}
}
\providecommand{\href}[2]{#2}

\end{document}